\documentstyle[amsmath,amssymb,amscd,graphicx,psfrag,epsf,amsthm]{article}

\newtheorem{thm}{Theorem}%[section]

\newtheorem{lem}[thm]{Lemma}
\newtheorem{cor}[thm]{Corollary}

\newtheorem{prop}[thm]{Proposition}

\newtheorem{conj}[thm]{Conjecture}
\theoremstyle{definition}
\newtheorem{defn}[thm]{Definition}

\newtheorem{say}[thm]{}
\newtheorem{exmp}[thm]{Example}

\newtheorem{rem}[thm]{Remark}          
\newtheorem{ack}{Acknowledgments}

\newtheorem{defn-thm}[thm]{Definition--Theorem}  %!!!!!!!!!!!!!!!!!!!!!!!!
\theoremstyle{remark}

\renewcommand{\c}[0]{{\mathbb C}}  

\renewcommand{\o}[0]{{\mathcal O}} 
\newcommand{\z}[0]{{\mathbb Z}}

\renewcommand{\r}[0]{{\mathbb R}} 

\renewcommand{\a}[0]{{\mathbb A}}

\newcommand{\p}[0]{{\mathbb P}}
\newcommand{\f}[0]{{\mathbb F}}

\newcommand{\map}[0]{\dasharrow}
\newcommand{\qtq}[1]{\quad\mbox{#1}\quad}
\newcommand{\spec}[0]{\operatorname{Spec}}

\newcommand{\br}[0]{\operatorname{Br}}
\newcommand{\var}[0]{\operatorname{Var}}
\newcommand{\hilb}[0]{\operatorname{Hilb}}

\def\into{\DOTSB\lhook\joinrel\rightarrow}

\begin{document}
\bibliographystyle{amsplain}

\title{Conics in the Grothendieck ring}
\author{J\'anos Koll\'ar}

\maketitle

The aim of this note is to describe the subring of the 
Grothendieck ring generated by smooth conics.
As a ring this is quite complicated, with many zero divisors,
but the
 description of the defining relations is entirely elementary.

\begin{defn} Let $k$ be a field. The {\it Grothendieck ring}
of $k$-varieties, denoted by $K_0[\var_k]$
is defined as follows.

Its additive group
 is the Abelian group whose generators are
the isomorphism classes of reduced, quasi projective $k$-schemes
and the relations are
$$
[X]=[Y]+[X\setminus Y],
$$
whenever $Y$ is a closed subscheme of $X$.

Multiplication is defined by $[X]\cdot [Y]=[X\times_k Y]$.
\end{defn}

The Grothendieck ring
of $k$-varieties is still very poorly understood. 
In characteristic zero, the quotient of 
$K_0[\var_k]$ by the ideal generated by $[\a^1]$ is naturally
isomorphic to the ring $\z[SB_k]$, where  $\z[SB_k]$ is the free abelian group 
 generated by   the stable birational
 equivalence classes
of smooth, projective, irreducible  $k$-varieties 
and multiplication is given by the product of varieties
\cite{la-lu}. (The cited  paper proves this
 over algebraically closed fields only, but the proof works
over any field of characteristic zero using the 
 birational factorization theorem as given in
 \cite[Remark 2 after Theorem 0.3.1]{akmw}.
Note also that the product of two irreducible  $k$-varieties 
is not necessarily irreducible, so  $\z[SB_k]$ is not a monoid ring
if $k$ is not algebraically closed.)

Zero divisors in the Grothendieck ring
of $\c$-varieties
were found by \cite{poonen}.

Here we  give further examples of nontrivial
behaviour of these rings by studying products of conics.
This gives interesting  examples only when the field $k$
is not algebraically closed.

\begin{thm}\label{main.thm} Let $k$ be a number field 
or the function field of an algebraic surface over $\c$. Let
 $C_i:i\in I$ and
$C'_j:j\in J$ be two collections of smooth   conics defined over $k$
(repetitions allowed).

Then $[\prod_{i\in I}C_i]= [\prod_{j\in J}C'_j]$ in the Grothendieck ring
iff $|I|=|J|$ and  the subgroup  of the Brauer group (cf.\ (\ref{brau.defn}))
generated by the first collection
$\langle C_i : i\in I\rangle \subset \br(k)_2$ 
 is the same as the  
subgroup  of the Brauer group 
generated by the second collection
$\langle C'_j : j\in J\rangle \subset \br(k)_2$.
\end{thm}

\begin{rem} (1)  The precise conditions on $k$ for the proof to work
are given in (\ref{field.cond}). These are satisfied for many other
fields, but fail for function fields of more than 2 variables.
It is not clear to me, however,  if any condition is needed on $k$.

(2)  Any isomorphism of two products 
$\prod_{i\in I}C_i$ and $\prod_{j\in J}C'_j$
is given in the obvious way: by a one--to--one map $g:I\to J$
and isomorphisms  $C_i\cong C'_{g(i)}$ for every $i$.

This can be proved many ways. Here is one using extremal
rays.

If $X$ is any projective variety, the cone of curves
of $X\times \p^1$ is generated by the cone of curves of
$X\cong X\times \{0\}$ and by $\{x\}\times \p^1$.
Using this repeatedly, we obtain that the
cone of curves of $(\p^1)^m$ is generated by the fibers of the
$m$ coordinate projections $(\p^1)^m\to (\p^1)^{m-1}$. Thus
the $|I|$ coordinate projections 
$$
\pi_{i'}: \prod_{i\in I}C_i\to \prod_{i\in I, i\neq i'}C_i 
$$
are ine one--to--one correspondence with the exremal rays of
$\prod_{i\in I}C_i$. Hence the  product structure can be recovered from
the intrinsic geometry of $\prod_{i\in I}C_i$.
\end{rem}

\begin{cor}  \label{gring.cor} Let $k$ be a  number field 
or the function field of an algebraic surface over $\c$.
The subring of the Grothendieck ring generated by smooth conics
is isomorphic to the ring generated by the isomorphism classes of
smooth conics modulo the ideal generated by the relations
$[C_1]\cdot [C_2]-[C_1]\cdot [C_1\ast C_2]$.

This ring can also be described as follows.

Let  $G\subset \br(k)_2$ be a finite subgroup with basis
$B_1,\dots,B_s$. Then $[B_1\times \cdots \times B_s]$
depends only on $G$ and it is denoted by $C(G)$.
The trivial subgroup gives $C(0)=[\spec k]$.

The Grothendieck ring of  conics is the free abelian group generated by
the elements $C(G)\cdot [\p^1]^m$ with multiplication
$$
C(G_1)\cdot C(G_2)=C(\langle G_1,G_2\rangle)\cdot 
[\p^1]^{\dim G_1+\dim G_2-\dim  \langle G_1,G_2\rangle}.
$$
where $\dim G$ denotes dimension as an $\f_2$ vector space.
\end{cor}

\begin{rem} The last description shows  that the 
Grothendieck ring of  conics
does not have nilpotents. Indeed, given an element 
$g=\sum \gamma_{G,m} C(G)\cdot [\p^1]^m$,
let $G_0$ be  minimal such that $\gamma_{G_0,m}\neq 0$ for some $m$, chosen
 also minimal.
Then the coefficient of 
$C(G_0)\cdot [\p^1]^{sm+(s-1)\dim G_0}$ in $g^s$ is 
$\gamma_{G_0,m}^s\neq 0$.
\end{rem}

The simplest example of nontrivial birational maps between
products of conics is the following. The whole description of the
Grothendieck ring of  conics is only 
a more elaborate version it.

\begin{exmp}\label{elementary.lem} Let $C$ be a smooth plain conic. Then
$C\times \p^1$ is birational to $C\times C$, and they have the same
class in the Grothendieck ring
of $k$-varieties.

Thus $[C]\cdot([\p^1]-[C])=0$ and $[C]$ is a zero divisor in
the  Grothendieck ring
of $k$-varieties if $C$ has no $k$-points.

Proof. $C\subset \p^2$ is a conic and we think of $\p^1$
as a line in the same $\p^2$.

Given $p,q\in C$, the line connecting them intersects
$\p^1$ in a point $\phi(p,q)$. 
$(p,q)\mapsto (p, \phi(p,q))$ gives a rational map
$C\times C\map C\times \p^1$. 
Conversely, given $p\in C$ and $r\in \p^1$ the line connecting them intersects
$C$ in a further point $\phi^{-1}(p,r)$.

Let $s,s'\in C(\bar k)$ be the two intersection points of $C$ and $\p^1$.
 $\phi$ is not defined at the pairs $(s,s')$ and $(s',s)$.
$\phi^{-1}$ is not defined at the pairs $(s,s)$ and $(s',s')$.
Easy computation shows that $\phi$ becomes an isomorphism
after we blow up the indeterminacy loci.
The blown-up surface is denoted by $B(C\times C)$.
As $k$-schemes, $(s,s')\cup(s',s)$ and  $(s,s)\cup(s',s')$ are
both isomorphic to $\spec_kk(s)$.

Thus $[C\times C]$ and $[ C\times \p^1]$ can both be written as
$$
[B(C\times C)]-[\p^1]\cdot [\spec_kk(s)]+[\spec_kk(s)].
$$

In order to see that  $[C]$ is a zero divisor in
the  Grothendieck ring
of $k$-varieties, we need to prove that
$[\p^1]-[C]$ is not zero.  
By \cite{la-lu}, it is sufficient to prove that 
$\p^1$ and $C$ are not stably birational. This is however easy, 
since having $k$-points is a stably birational invariant.
\end{exmp}

\begin{say}[Products of conics and the Brauer group]{\ }
\label{brau.defn}

Below we  give an elementary geometric description
of the Brauer group of conics, denoted by $\br(k)_2$.
(For the fields that we are considering, this is the 2--torison subgroup
of the  Brauer group $\br(k)$. See \cite[X.4--7]{serre-lf} for a good
introduction and basic properties.)

Let $k$ be a field and $C_1,C_2$ 
two smooth conics defined over $k$.
The Brauer product of the two conics is defined as follows.
(I warn the reader in advance that
this definition only works because on a conic
the Hilbert scheme of points is isomorphic to the Hilbert scheme
of degee 1 divisors. On a higher dimensional Severi--Brauer variety
the Hilbert scheme of points is  dual  to the Hilbert scheme
of degee 1 divisors and one has to distinguish these systematically.)

Start with $C_1\times C_2$. As a first approximation, we
construct a 3--dimensional variety, denoted by   $P(C_1,C_2)$.
We would like to say that   $P(C_1,C_2)$ is
the 3--dimensional ``linear system'' of divisors of bidegree $(1,1)$
on $C_1\times C_2$.
The problem is that in general no such divisor is
defined over $k$. Thus we look at the linear system
$|-K|$ where $K=K_{C_1\times C_2}$ is the canonical class.
This corresponds to divisors of bidegree $(2,2)$. Then
$ P(C_1,C_2)\subset |-K|$ is the subscheme consisting of
those divisors which are everywhere double. Over
the algebraic closure of $k$ we recognize this as the (doubled)
elements of the linear system $|\o(1,1)|$. 

Alternatively, the Hilbert scheme $\hilb(C_1\times C_2)$
has an irreducible component parametrizing  divisors of bidegree $(1,1)$.
This is again $P(C_1,C_2)$.

Thus $P(C_1,C_2)$ is a 3--dimensional $k$-variety which is
isomorphic to $\p^3$ over $\bar k$.

There is a natural embedding $C_1\times C_2\into  P(C_1,C_2)$
where we map a point $(p,q)\in C_1\times C_2$
to the divisor $2(\{p\}\times C_2+C_1\times \{q\})$.

In general this is all one can do. Ther are, however, important cases
when such a product  $P(C_1,C_2)$ contains a degree 1 smooth
curve (a line over $\bar k$) defined over $k$. In this case I call this 
degree 1 curve the Brauer product
of $C_1$ and $C_2$ and   denoted it by $C_1\ast C_2$.
(The terminology ``Brauer product'' does not seem to be standard.)

It turns out that this is well defined.

To see this,
let $P$ be a  3--dimensional $k$-variety which is
isomorphic to $\p^3$ over $\bar k$. Let $L_1,L_2\subset P$
be degree 1 smooth curves defined over $k$ 
and let $L'\subset P$ be another such curve
disjoint from both.
(Over an infinite field we can obtain $L'$ as the image of $L_1$
by a general automorphism of $P$.)
 Then $L_1$ and $L_2$ are both isomorphic to the
Hilbert scheme of degree 1 surfaces containing $L'$.

$C\ast C$ is always defined and it is isomorphic to $\p^1$.
Indeed, the diagonal  $\Delta\subset C\times C$ is defined over $k$
thus $P(C,C)$ is $k$-isomorphic to $\p^3_k$. Hence the Brauer group of
conics is a 2--group.
\end{say}

\begin{lem}\label{field.cond}
For a field $k$ the following two conditions are equivalent.
\begin{enumerate}
\item The  Brauer product of 2 smooth conics is again a conic.
\item For any two smooth conics $C_1,C_2$ defined over $k$ there is a
degree 2 extension $k'/k$ such that both $C_1$ an $C_2$ have
$k'$-points.
\end{enumerate}
\end{lem}

Proof.  Let $L\subset P(C_1,C_2)$ be a degree 1 curve defined over $k$.
Then $L\cap (C_1\times C_2)$ is a degree 2 subscheme defined over $k$
with residue field $k'$.
By projection to the factors, $C_1,C_2$ both have points in
$k'$.

Conversely, if $C_1,C_2$ both have points in
$k'$ then so does their product. The unique line
in $P(C_1,C_2)$ passing through a $k'$ point is defined
over $k$. 
\qed
\medskip

The following result is well known in various forms,
see for instance \cite[p.209]{artin} or \cite[Thm.5.7]{sarkisov}

\begin{prop}\label{field.prop}
 The conditions in (\ref{field.cond})
hold in  the following two cases:
\begin{enumerate}
\item $k$ is a number field.
\item $k$ is the function field of an algebraic surface over an algebraically 
closed field. More generally, for $C_2$-fields.
\end{enumerate}
\end{prop}

Proof. Here is a geometric version of some of the classical proof.

Let $G(1,P(C_1,C_2))$ denote the Grassmannian of lines in $P(C_1,C_2)$.
 We need to prove
that it has a $k$-point.

More generally,
let $P$ be a $k$-variety which is isomorphic to $\p^n$ over $\bar k$
and assume that there is a quadric hypersurface  $Q\subset P$ defined over $k$.
As explained in \cite[4.5]{artin} 
the Grassmannian of lines  $G(1,P)$ is embedded
into $\p^{\binom{n+1}{2}-1}$ the usual way.

For $n=2$ the  Grassmannian of lines  $G(1,P)$ is thus a quadric
in $\p^5$, and so it has a point over any $C_2$ field, proving the second
part.

If $k=\r$ then  a $\c$-point of $P$ and its conjugate determine
a real line, so $G(1,P)(\r)\neq\emptyset$. 
Thus $G(1,P)$ is a quadric in 6 variables which has a point
in all real completions of $k$.
Therefore  $G(1,P)$ has a $k$-point by the Hasse--Minkowski theorem.
\qed

\begin{say}[Proof of (\ref{main.thm})]{\ }

The key point  to show one direction 
is the following generalization of (\ref{elementary.lem}): 

\begin{lem}\label{2conics.lem}
 Let $C_1,C_2$ be smooth conics
such that their Brauer product $C_1\ast C_2$ is again a  smooth conic. Then
\begin{enumerate}
\item $C_1\times C_2$ is birational to $C_1\times (C_1\ast C_2)$.
\item $[C_1\times C_2]=[C_1\times (C_1\ast C_2)]$ in the Grothendieck ring.
\end{enumerate}
\end{lem}

Proof. First we write down a rational map
$\phi:C_1\times C_2\map C_1\ast C_2$. Then we check that
$\phi$ and the fist projection $\pi_1:C_1\times C_2\to C_1$
give a birational map
$$
(\pi_1,\phi):C_1\times C_2\map C_1\times (C_1\ast C_2).
$$
Finally we see that this gives an identity in the Grothendieck ring.

{\it Geometric description.} By assumption there is a degree 2 point
$Q\in C_1\times C_2\subset P(C_1,C_2)$.
Let $L'$ be the unique degree 1 curve through $Q$ and let
$L\in P(C_1,C_2)$ be any degree 1 curve  disjoint
from $L'$. Projection from $L'$ to $L$ gives $\phi$.

{\it Algebraic description.} Assume for simplicity that the
characteristic is different from 2.
If the common point is in the field
$k(\sqrt{a})$, we can assume that the conics are given by
equations
$$
C_1=(x_1^2-ax_2^2-bx_3^2=0)
\qtq{and}
C_2=(y_1^2-ay_2^2-cy_3^2=0).
$$
Then their Brauer product can be given as
$$
C_1\ast C_2=(z_1^2-az_2^2-bcz_3^2=0)
$$
and $\phi$ is given by
$$
(z_1:z_2:z_3)=(x_1y_1+ax_2y_2: x_1y_2+x_2y_1: x_3y_3).
$$

$\phi^{-1}$ is obtained as follows. Pick a point $p\in C_1$ and $r\in L\cong
C_1\ast C_2$. $\{p\}\times C_2$ embeds as a line into $P(C_1,C_2)$
and $\phi^{-1}(p,r)$ is the intersection point of
this line with the plane $\langle L',r\rangle$ spanned by $L'$ and $r$. 
This is not defined only
if $\{p\}\times C_2\subset \langle L',r\rangle$.
This happens exactly when $\langle L',r\rangle$ is one of the two
tangent planes of $C_1\times C_2$ at a point of $Q$
and $\{p\}\times C_2$ is the corresponding line through
that point of $Q$.

Thus we see that 
 $C_1\times C_2$ becomes isomorphic  to $C_1\times (C_1\ast C_2)$
after we blow up subschemes isomorphic to $Q$ in both of them.
As in (\ref{elementary.lem}) this shows that
 $[C_1\times C_2]=[C_1\times (C_1\ast C_2)]$.
\qed
\medskip

Assume that the subgroup  $G_I\subset \br(k)_2$ generated by
the $C_i$-s is the same as the  
subgroup $G_J\subset\br(k)_2$ generated by
the $C'_j$-s.
Fix a minimal generating set $\{B_s:s\in S\}$ of $G$.
By a simple group theoretic lemma (\ref{group.lem}) and a repeated
 application of
(\ref{2conics.lem}),
 $\prod_{i\in I}C_i$ is birational to
$$
 (\p^1_k)^{|I|-|S|}\times \prod_{s\in S}B_s,
$$
and they have the same class in the Grothendieck ring.
The same holds for  $\prod_{j\in J}C'_j$.
Thus $\prod_{i\in I}C_i$  and $\prod_{j\in J}C'_j$ are birational
and they have the same class in the Grothendieck ring.

Conversely, assume that $G_I\neq G_J$. We may assume that
$G_J\not\subset G_I$ and so there is an index $j_0$ such that
the class of $C'_{j_0}$ is not in $G_I$. We claim that in this case there is
no rational map from $\prod_{i\in I}C_i$  to $C'_{j_0}$, hence
no rational map from $\prod_{i\in I}C_i$  to $\prod_{j\in J}C'_j$.
Thus they are not birational and not even stably birational,
hence they represent different elements of the
Grothendieck ring by \cite{la-lu}.

The proof is by induction on $|I|$, the case $|I|=0$  being clear.
Pick $i_0\in I$ and set $I':=I\setminus \{i_0\}$ and $K=k(C_{i_0})$.
By (\ref{amits.lem}),  the kernel of
$G_I\to \br(k)_2$ is generated by $C_{i_0}$ and so 
 the class of $C'_{j_0}$ in 
$\br(k)_2$
is not in the subgroup $G_{I'}\subset \br(K)_2$ generated by
the $C_i$-s for $i\in I'$. 
By induction, there is no $k(C_{i_0})$-map 
from $\prod_{i\in I'}C_i$  to $C'_{j_0}$,
and so no $k$-map from $\prod_{i\in I}C_i$  to $C'_{j_0}$.
This completes the proof of (\ref{main.thm}).

In the proof we used only the relations given by
(\ref{2conics.lem}), and this gives the first description
of the Grothendieck ring in (\ref{gring.cor}).

It is clear the all the possible products
 $\prod_{i\in I}C_i$ 
generate the Grothendieck ring as an additive abelian group, and we
have established that each such product is identical to a unique element
of the form
$$
 [\p^1_k]^{|I|-|S|}\cdot \prod_{s\in S}[B_s]=[\p^1_k]^{|I|-|S|}\cdot C(G).
$$
This gives the second description in (\ref{gring.cor}).
\qed

\begin{lem}\label{group.lem} Let $G$ be a finite abelian 2--group
with a minimal generating set $b_1,\dots,b_m$.

Let $e_1,\dots,e_s$ be any generating collection of elements of $G$,
 repetitions allowed.
Then $e_1,\dots,e_s$ can be transformed into the 
collection $b_1,\dots,b_m, 0,\dots,0$ by repeated application of the
following operation:
$$
\mbox{Pick $e_i,e_j$ and replace $e_j$ by $e_i+e_j$.}\qed
$$
\end{lem}

The following is a very special case of an old result of \cite{amitsur}.

\begin{lem}\label{amits.lem}
 Let $C,C'$ be smooth conics and
$g:C\to C'$ a rational map. Then either $C'\cong \p^1$
or $C\cong C'$.

Therefore, if $C'\not\cong C,\p^1$ then $C'$ does not have
a  $k(C)$-point.
\end{lem}

Proof. Let $G\subset C\times C'$ be the graph of $g$.
It is a divisor of bidegree $(1,\deg g)$. A class of bidegree
$(0,2)$ is defined over $k$, so we obtain that either
the linear system $|\o(1,0)|$ or the linear system $|\o(1,1)|$
has a member over $k$. In the first case $C'$ has a $k$-point and
$C'\cong \p^1$ and in the second case we get a graph of an isomorphism.\qed
\end{say}

The birational part of (\ref{main.thm}) is easy to
state and prove for higher dimensional Severi--Brauer
varieties:

\begin{prop} Two products $\prod_{i\in I}P_i$ and $\prod_{j\in J}P'_j$ 
of Severi--Brauer
varieties are stably birational iff  the subgroup 
 $\langle P_i : i\in I\rangle \subset \br(k)$ 
 is the same as the  
subgroup  $\langle P'_j : j\in J\rangle \subset \br(k)$.
\end{prop}

It is less clear how to formulate the
Grothendieck ring version in general.  Even for products of
2--dimensional Severi--Brauer
varieties I found it difficult to write down a suitable analog
of (\ref{2conics.lem}). Additional problems arise
when the two Severi--Brauer
varieties do not have a common splitting field.

\begin{ack}  
Partial financial support was provided by  the NSF under grant number 
DMS02-00883. 
\end{ack}

\vskip1cm

\noindent Princeton University, Princeton NJ 08544-1000

\begin{verbatim}kollar@math.princeton.edu\end{verbatim}


\begin{thebibliography}{Chislenko88}
\bibitem[Poonen]{poonen} B.\ Poonen, 
 The Grothendieck ring of varieties is not a domain. math.AG/0204306
\end{thebibliography}
\end{document}